\documentstyle[12pt]{article}
\def\be{\begin{equation}}
\def\ee{\end{equation}}
\def\beq{\begin{equation}}
\def\eeq{\end{equation}}
\def\bd{b^{\dagger}}
\def\bda{b^{\dagger}_{\alpha}}

\def\bb{b_{\beta}}
\def\bdm{b^{\dagger}_{\mu}}

\begin{document}

\begin{center}
{\bf QUONS AS $su(2)$ IRREDUCIBLE TENSOR OPERATORS}
\end{center}

\begin{center}
{\it S.S. Avancini,F.F.de Souza Cruz,
J.R.Marinelli and D.P.Menezes\\
{\it Depto de F\'{\i}sica - CFM - Universidade Federal de Santa Catarina - Florian\'opolis - SC - CP. 476 - CEP 88.040 - 900 - Brazil}}
\end{center}

\vspace{0.50cm}

\begin{abstract}
We prove that, for the quon algebra, which interpolates between the Bose and Fermi 
statistics and depends on a free parameter $q$,it is possible to build
an $su(2)$ irreducible representation.
One of the consequences of this fact
is that the quons couple via the same angular momentum coupling rules
obeyed by ordinary bosons and fermions.

\end{abstract}

\newpage

\section {Introduction}

The {\it quons} are particles whose statistics interpolates
between the boson and fermion statistics and depend on a special deformation
parameter which varies from +1 to -1 \cite{green}.The usual bosonic 
(fermionic) algebra is recovered when the deformation parameter is equal 
to $ +1 (-1)$. 
The $q-$deformed
commutation relation obeyed by the {\it quons} comes out as an specific
case of
the commutation relations obeyed by the $q-$deformed boson operators
\cite{qa} of the
quantum algebras \cite{kibler} when we consider just one degree of freedom
for the system. In general, however, they lead to different commutation
relations and very specific properties, as discussed, for example, in
reference \cite{chaichian}.
The quon algebra comes out\cite{green} naturally in the study of the 
interactions of particles having small violations of Fermi and Bose 
statistics. Limits for these violatios\cite{green,newgreen} have been
obtained in atomic and nuclear physics experiments.

In addition to these speculative works, the quon algebra has been employed
as a tool to describe many-body problems involving composite 
particles \cite{avga}. 
In these problems, due to the simultaneous presence of composite 
and constituinte degrees of freedom, the creation and annihilation operators 
of the composite particles deviates from the canonical ones.
 The deformation parameter $q$ of the quon algebra may be considered as a 
measure of the effects of the statistics of the internal degrees of 
freedom of the composite particles.

Besides the intrinsic interest in the study of {\it quons} and of the
fundamental behaviour of systems of particles that obey its statistics
\cite {Newgreen}, the possibility of describing particles with properties
which interpolate
between fermion and boson ones is quite appealing in several areas of
physics, specially in the treatment of many body systems through boson
expansions
\cite{nosso1},\cite{nosso2}.
In the above mentioned works, as in many other possible applications,
we have to deal with tensor operators that have to be coupled in order to
define the physical operators of interest. In quantum algebras, it has been
shown that q-deformed Clebsh-Gordan coefficientes must be introduced,
which are now available in the literature (see \cite{qCG}, for example).
However they are not suitable for {\it quons}, once as we said before,
the commutation relations and therefore the algebraic structure obeyed by 
quons are different in general from the ones that define the quantum algebras.

In the present work we construct {\it quon} operators which behave as usual
$su(2)$ irreducible tensor operators. For this purpose we find
a realization of the $su(2)$ algebra in terms of the $q$--deformed
bosons by building the su(2) generators $J_+$, $J_-$ and $J_z$
\cite{livroBie} as functions of the quon operators. Within the
quon algebra the number operators are written as series in terms of
the deformation parameters. The $su(2)$ generators built in this way also
commute with the operators used to construct the hamiltonian.

\vspace{0.5cm}

\section{The quon algebra}

We assume that $b_m$,$b^{\dagger}_m$, $m=-j,...+j$ are (2j+1) operators which 
satisfy the $q$--mutation relations (or the quon algebra) defined by
\cite{green}, i.e.,
 
\be
[b_m,b^{\dagger}_{m^{\prime}}]_{q}=
b_m ~{\bd}_{m^{\prime}}-q ~{\bd}_{m^{\prime}} ~b_m = 
 ~{\delta}_{m,m^{\prime}}.
\label{comu}
\ee
This relation is a deformation of the Bose and Fermi algebras and
interpolates between those algebras when $q$ goes from +1 to -1.
The Fock space is constructed as usual from the application of the 
creation operators, ($b_\alpha^\dagger$), on the vacuum state, 
defined as $b_\alpha|0>=0$ for all $\alpha$. It has been 
shown\cite{fivel,zagier} that the squared norm of any polynomials of the 
creation operators, $b_\alpha^\dagger$,
acting on the vacuum state is positive definite.

In what follows we derive some important results in order to construct  
the $su(2)$ generators. We have at our disposal
the transition number operator \cite{green} $N_{\alpha\beta}$, 
which has the usual commutation relations:
\be
[N_{\alpha\beta},\bdm]={\delta}_{\beta\mu}~\bd_{\alpha}~,~
[N_{\alpha\beta},b_\mu]=-{\delta}_{\alpha\mu}~b_{\beta}~~. \label{tran}
\ee
The transition number operator is an infinite series in the $\bdm$'s
and $b_\mu$'s and it has been obtained in closed form in 
\cite{zagier,mel}. The general structure of this operator is \cite{mel} :
\be
N_{\alpha\beta}~=~\bda \bb~+~\sum_{n=1}^{\infty}~\sum_{(i_1 ...i_n)}~
\sum_{\pi}~c_{\pi(i_1),...,\pi(i_n), i_1 ...i_n} ~ 
(Y_{\alpha\pi(i_1), ...\pi(i_n)})^{\dagger} ~ Y_{\beta i_1 ...i_n}
\label{nger}
\ee
where the summation over $\pi$ means that we consider all different 
permutations of the indices $i_1,i_2, ...,i_n$, including their repetitions.
$\pi(i_1),\pi(i_2), ...\pi(i_n)$ corresponds to each of these permutations
 and $Y^\dagger$ is the adjoint of $Y$.
The operators $Y_{\beta i_1 ...i_n}$ are obtained through the recursion 
relations \cite{mel,dori}:
\be
Y_{k i_1 ...i_{n+1}}~=~Y_{k i_1 ...i_n}~b_{i_{n+1}}~-~
q^{n+1}~b_{i_{n+1}}~Y_{k i_1 ...i_n}
\label{recur}
\ee
with
\[
Y_{ki}~=~b_k b_{i}~-~q~b_i b_k
\]
The coefficients $c_{\pi(i_1),...,\pi(i_n), i_1 ...i_n}$ 
are only functions of the deformation parameter $q$ \cite{mel,dori}.
The first terms in the series are:
\be
N_{\alpha \beta} = \bd_{\alpha} b_{\beta} + (1-{q}^2)^{-1}
\sum_m(\bd_m\bd_{\alpha} - q \bd_{\alpha}\bd_m)(b_{\beta}b_m-q b_m b_{\beta}) 
+... ~~.
\ee
From the recursion relation, eq.(\ref{recur}), and eq.(\ref{tran}),
one can show that
\[
[N_{\alpha \beta}~,(Y_{{\alpha}^{\prime} i_1 ...i_n})^{\dagger} ]=
\delta_{\beta i_1} (Y_{{\alpha}^{\prime} \alpha i_2 ...i_n })^{\dagger}+
\delta_{\beta i_2} (Y_{{\alpha}^{\prime} i_1 \alpha i_3 ...i_n})^{\dagger} 
+...+
\delta_{\beta i_{n-1}} (Y_{{\alpha}^{\prime} i_1 i_2 ...\alpha i_n})^{\dagger}
\]
\be
+~ \delta_{\beta i_n} 
(Y_{{\alpha}^{\prime} i_1 i_2 ...i_{n-1}\alpha})^{\dagger}~+~
\delta_{\beta \alpha^\prime } (Y_{ \alpha i_1 i_2 ...i_n} )^{\dagger}.
\label{ncom}
\ee
As a subsidiary result, it follows from eq.(\ref{ncom})  that;
\[
[N_{\alpha \beta}~,(Y_{\alpha^\prime \pi(i_1) ...\pi(i_n) })^{\dagger}
Y_{\beta^\prime i_1 ... i_n } ] =
\delta_{\beta \alpha^\prime }~
(Y_{\alpha \pi(i_1) ...\pi(i_n) })^{\dagger}
~ Y_{ \beta^\prime i_1 ... i_n } 
\]
\be ~-~
\delta_{\alpha \beta^\prime }~
(Y_{ \alpha^\prime \pi{i_1} ...\pi{i_{n+1}}})^{\dagger}
~ Y_{\beta i_1 ... i_{n+1} } \label{NYY}.
\ee

In what follows we calculate the commutator
$[N_{\alpha \beta},N_{\alpha^{\prime} \beta^{\prime}}]$. For this
purpose, we substitute
$N_{\alpha^{\prime} \beta^{\prime} }$ given in eq.(\ref{nger})
into the commutator above, which yields
\[
[N_{\alpha \beta}, N_{\alpha^{\prime} \beta^{\prime}}] =
[N_{\alpha \beta},
~{b^{\dagger}}_{\alpha^{\prime}} ~b_{\beta^{\prime}}]~
\]
\[
+~\sum_{n=1}^{\infty}~\sum_{(i_1 ...i_n)}~
\sum_{\pi}~c_{{\alpha^{\prime}}\pi(i_1),...,\pi(i_n),
{\beta^{\prime}} i_1 ...i_n} ~ 
[N_{\alpha ,\beta},(Y_{{\alpha^{\prime}} \pi(i_1), ...\pi(i_n)})^{\dagger} ~ 
Y_{{\beta^{\prime}} i_1 ...i_n}].
\]
Using eqs.(\ref{tran}) and (\ref{NYY}), we obtain an important 
relation obeyed by the transition number operator:
\be
[N_{\alpha \beta},N_{\alpha^{\prime} \beta^{\prime}}] =
  {\delta}_{\beta {\alpha}^{\prime}}~N_{\alpha\beta^{\prime}}
-{\delta}_{\alpha {\beta}^{\prime}} N_{{\alpha}^{\prime}\beta } ~~.
\label{rip} 
\ee
The above commutation relation shows that the 
transition number operators are the generators of an $su(2j+1)$ algebra. 

\section{ $su(2)$ tensor operators}
In order to see that the $\bdm$ is an irreducible tensor 
operator, we have to show that 
\[
[J_{0} ,\bdm ]= \mu~\bdm   ~~, 
\]
\be
[J_{\pm},\bdm ]= \sqrt{(j\mp \mu)(j\pm \mu +1) } \bd_{\mu \pm 1}~~,\label{tens}
\ee
where $J_{0} ,J_{\pm}$ are usual $su(2)$ generators, 
\[
[J_{0},J_{\pm}]= {\pm} J_{\pm}  ~~, 
\]
\be
[J_{+} , J_{-} ]= 2 J_{0} ~~,
\ee 
which have to be built from the same $\bdm$ operators.
Notice that the label $m$ corresponds to the projection of
$j$ and we are restricted to a single $j$--level. Operators
belonging to different levels commute with each other.
Thus we may define the $su(2)$ generators as:
\[
J_0 = \sum_{\nu=-j}^{+j} \nu~N_{\nu \nu}  
\]
\be
J_{\pm} = \sum_{\nu=-j}^{+j} \sqrt{(j\mp\nu)(j\pm\nu+1)} N_{\nu\pm 1~\nu} \label{ma}
\ee
where the index $\nu$ takes the values $-j, -j+1,...j-1,j$ .
It is then simple to verify that the $\bdm$ operators
satisfy eqs.(\ref{tens}) and then form a genuine irreducible  
tensor operator of the  $su(2)$ algebra:
\[
[J_0,\bdm]= \sum_{\nu} \nu [N_{\nu \nu},\bdm]= \sum_{\nu}{\delta}_{\nu\mu}
\bdm = \mu \bdm
\]
and
\[
[J_{\pm},\bdm]= \sum_{\nu} \sqrt{(j\mp\nu)(j\pm\nu+1)} [N_{\nu\pm 1~\nu},\bdm] =
\sum_{\nu} \sqrt{(j\mp\nu)(j\pm\nu+1)} {\delta}_{\nu\mu}{\bd}_{\nu\pm 1} = 
\]
\[ 
\sqrt{(j\mp\mu)(j\pm\mu+1)} {\bd}_{\mu\pm 1} ~~.
\]
To finish our demonstration we still have to show that the 
operators defined in eqs.(\ref{ma}), are indeed 
generators of the usual $su(2)$ algebra. This proves the 
consistency of our approach. From eq.(\ref{ma}), we can write
\be
[J_0,J_{\pm}]= \sum_{\nu,{\nu}^{\prime}}{\nu} 
\sqrt{(j\mp{\nu}^{\prime})(j\pm{\nu}^{\prime}+1)}
[N_{\nu \nu} ,N_{{\nu}^{\prime}\pm 1 ~{\nu}^{\prime} } ]
\ee
\[ = \sum_{\nu~{\nu}^{\prime}} 
\nu~ \sqrt{(j\mp {\nu}^{\prime})(j\pm {\nu}^{\prime}+1)}
({\delta}_{\nu~ {\nu}^{\prime}\pm 1} N_{\nu{\nu}^{\prime}} 
-{\delta}_{\nu~{\nu}^{\prime}} N_{{\nu}^{\prime}\pm 1~\nu } ) = \pm J_\pm ~~.
\]

One can show in an analogous way that $[J_+ , J_- ]=2~J_0 $
and therefore the usual $su(2)$ generators are consistent 
with the quon algebra, eq.(\ref{comu}).
\section{ Concluding Remarks} 

Thus, the tensor coupling 
has to be done with the usual Clebsch-Gordan coefficients.

In the Introduction we have mentioned, as an application of the
{\it quon} statistics, the description of many body states by means of
bosons, as it is done, for example in the IBM \cite {ibm}. In that model
the so called $s$ and $d$ bosons are introduced
and in more sophisticated versions higher angular momentum bosons
are also needed. It follows from the above demonstration that we
can now define deformed $s$ and $d$ bosons, which behave like {\it quons}
but keep the same angular momentum coupling rules and coefficients as
in the non-deformed case. Of course, many other applications are
possible.

\vskip 0.35in
 This work has been partially supported by CNPq.

\end{document}